\documentclass[11pt,openbib]{article}
\usepackage{amsmath,amssymb}

\newcommand{\R}{\mathbb{R}}

\newcommand{\spp}{\mathop{\mathrm{sp}}\nolimits}

\newcommand{\Spp}{\mathop{\mathrm{Sp}}\nolimits}

\newcommand{\Ad }{\mathop{\mathrm{Ad}}\nolimits}

\newcommand{\trace}{\mathop{\mathrm{tr}}\nolimits}
\newcommand{\setrule}{\, \rule[-4pt]{.5pt}{13pt}\, }

\newcommand{\rowspace}{\rule{0pt}{16pt}}

\newcommand{\spann}{\mathop{\rm span}\nolimits}

\allowdisplaybreaks

\begin{document}
\begin{center}
{\Large \bf Coadjoint orbits of the odd real symplectic group}\\
\rule{0pt}{20pt} Richard Cushman\footnotemark 
\end{center} \bigskip 
%\addtocounter{footnote}{1}
\footnotetext{Department of Mathematics and Statistics, 
University of Calgary, Calgary, Alberta, T2N 1N4 Canada. email: r.h.cushman@gmail.com}
\addtocounter{footnote}{1}
\footnotetext{version: \today}
\vspace{-.3in}\begin{abstract}
We give a representative of every coadjoint orbit of the odd symplectic group. Our argument 
follows that used for the Poincar\`{e} group but the details differ.
\end{abstract}

 \medskip 

 \bigskip 

We introduce the notion of cotype and reduce the
classification of coadjoint orbits of the odd real symplectic group to 
the classification of types for the Lie algebra of the real symplectic group, see \cite[p.352-3]{burgoyne-cushman}. For a classification of the adjoint orbits of the odd symplectic group see 
\cite{cushman}. We follow the line of argument in \cite{cushman-vanderkallen} for classifying 
the coadjoint orbits of the Poincar\`{e} group, but there are some new features. 

\section{The odd symplectic group}
%%%%%%%%%%%%%%

We begin by defining some basic concepts. Let $(V,\omega )$ be a real \emph{symplectic vector space}. 
A real linear map $P:V \rightarrow V$, which preserves the symplectic form $\omega $, that is, 
$P^{\ast }\omega = \omega $, is called a linear \emph{symplectic map}. The set of all real 
linear symplectic maps is the Lie group $\Spp (V, \omega )$. Its Lie algebra $\spp (V, \omega )$ 
is the set of real linear maps $Y: V \rightarrow V$ such that $\omega (Yv,w) + \omega (v,Yw) =0$  
for every $v$, $w \in V$. Let $v^{\circ}$ be a nonzero vector in $V$. Let ${\Spp (V, \omega )}_{v^{\circ}}$ 
be the Lie group of all real linear symplectic maps of $(V, \omega )$ into itself, which leave 
$v^{\circ}$ fixed. ${\Spp (V, \omega )}_{v^{\circ}}$ is called the \emph{odd symplectic} group, 
see \cite{bates-cushman}. Its Lie algebra ${\spp (V, \omega )}_{v^{\circ}}$ consists of all 
$Y \in \spp (V, \omega )$ such that 
$Yv^{\circ} =0$. Let $Y = S+N$ be the Jordan decomposition of $Y$ into 
a sum of a semisimple linear map $S$ and a commuting nilpotent linear map $N$, which both lie in 
$\spp (V, \omega )$. Because $S$ and $N$ are polynomials in $Y$ with real coefficients and no constant terms. Since $Yv^{\circ} = 0$, it follows that $Sv^{\circ} = 0 = Nv^{\circ}$. So $S$ and $N$ lie in 
${\spp (V, \omega )}_{v^{\circ}}$. 

\section{Classification of coadjoint orbits}
%%%%%%%%%%%%%%%%%%%%%

A \emph{tuple} $(V,Y,v;\omega )$ is a real symplectic vector 
space $(V,\omega )$, a vector $v\in V$, and a real linear 
map $Y \in \spp (V, \omega )$. Two tuples $(V,Y,v; \omega)$ 
and $(V',Y',v'; {\omega }')$ are \emph{equivalent} if there 
is a bijective real linear map $P:V \rightarrow V'$ such that 
i) $Pv = v'$, ii) $P^{\ast }{\omega }' = \omega $, and iii) there 
is a vector $w \in V$ so that $Y' =P(Y + L_{v,w})P^{-1}$. Here 
$L_{v,w} = v \otimes w^{\ast } + w \otimes v^{\ast }$, where 
$v^{\ast }$ is the linear function on $V$ which 
sends $u$ to $\omega (v,u)$. Note that 
$L_{v,w} = L_{w,v}$ for every $v,w \in V$. \medskip 

\noindent \textbf{Fact 1.} For every $v,w \in V$, we have $L_{v,w} \in 
\spp (V,\omega )$. \medskip 

\noindent \textbf{Proof.} For $x,y \in V$ 
\begin{align}
\hspace{-.15in}\omega \big( (v \otimes w^{\ast })x + (w \otimes v^{\ast })x, y \big) 
+ \omega \big( x, (v \otimes w^{\ast })y + (w \otimes v^{\ast })y \big) 
& = \notag \\
& \hspace{-3.5in} = \omega \big( \omega (w,x)v + \omega (v,x)w, y \big) + 
\omega \big( x, \omega (w,y)v + \omega (v,y)w \big) \notag \\
& \hspace{-3.5in} = \omega (w,x) \omega (v,y) + \omega (v,x) \omega (w,y) 
+ \omega (x,w) \omega (v,y) + \omega (x,v) \omega (w,y) \notag \\
& \hspace{-3.5in} = 0, \quad \mbox{since $\omega $ is skew symmetric.} \notag 
\end{align}
Therefore $L_{v,w} \in \spp (V, \omega )$. \hfill $\square $ \medskip  

\noindent \textbf{Fact 2.} If $P\in \Spp (V, \omega )$, then $PL_{v,w}P^{-1} 
= L_{Pv,Pw}$. \medskip 

\noindent \textbf{Proof.} For $z \in V$ we have 
\begin{align}
(PL_{v,w}P^{-1})z & = PL_{v,w}(P^{-1}z)  
= P\big( (v \otimes w^{\ast })P^{-1}z + 
(w \otimes v^{\ast })P^{-1}z \big) \notag \\
& = w^{\ast }(P^{-1}z)Pv + v^{\ast }(P^{-1}z)Pw \notag \\
& = \omega (w, P^{-1}z)Pv + \omega (v,P^{-1}z)Pw \notag \\
& = \omega (Pw,z)Pv + \omega (Pv,z)Pw, \quad \mbox{since 
$P \in \Spp (V,\omega )$} \notag \\
& = (Pv \otimes (Pw)^{\ast })z + (Pw \otimes (Pv)^{\ast})z = 
L_{Pv,Pw}(z). \tag*{$\square $}
\end{align}  

Being equivalent is an equivalence relation on the collection 
of all tuples. An equivalence class of tuples is called a 
\emph{cotype}, which is denoted by $\nabla $. Suppose that 
the tuple $(V,Y,v; \omega )$ represents the cotype $\nabla $. 
If $v \ne 0$, then $\nabla $ is an \emph{affine} cotype. 
Define the \emph{dimension} of $\nabla $ to be $\dim V$. The \emph{height} of the 
cotype $\nabla$ is the largest positive integer $n$ such that $N^{n+1}V = \{ 0 \}$, 
where $N$ is the nilpotent summand in the Jordan decomposition of $Y$.

\medskip 

\noindent \textbf{Lemma 3.} Every affine cotype $\nabla $ of dimension 
$2n+2$ has a representative tuple $({\R }^{2n+2}, Y', e_{2n+1}; J)$. \medskip 

\noindent \textbf{Proof.} Let $(V,Y,v;\omega )$ be a tuple 
which represents the cotype $\nabla $. Let $\mathfrak{f} = 
\{ v_0,\, v_1, \ldots , v_{2n+1} = v \}$ be a basis of $V$ with 
respect to which the matrix of $\omega $ is $J =${\tiny 
$ \left( \begin{array}{rcc} 0 & 0 & 1 \\ 0 & \widetilde{J} & 0 \\
-1 & 0 & 0 \end{array} \right) ,$} where $\widetilde{J} =${\footnotesize 
$\left( \begin{array}{rc} 0 & I_n \\ -I_n & 0 \end{array} \right) $}. 
Let $P:V \rightarrow {\R }^{2n+2}$ be the bijective real linear 
map which sends the basis vector $v_i$ of $V$  to the standard basis vector $e_i$ of ${\R }^{2n+2}$ for $i =0, \ldots , 2n+1$. Then 
$P^{\ast }J = \omega $ and $Pv = e_{2n+1}$. The matrix $Y'$ of $Y$ with respect to the standard basis 
$\mathfrak{e} = \{ e_0, \, e_1, \ldots , e_{2n}, \, e_{2n+1} \} $ 
of ${\R }^{2n+2}$ is 
\begin{displaymath}
Y' = PYP^{-1} = P(Y +L_{0,v})P^{-1}. 
\end{displaymath}
Thus the tuple $({\R }^{2n+2}, Y', e_{2n+1}; J)$ is equivalent 
to the tuple $(V,Y,v;\omega )$. \hfill $\square $ \medskip 

\vspace{-.15in}\noindent \textbf{Proposition 4.} For $Y \in \spp ({\R }^{2n+2}, J)$ let 
${\ell }_Y$ be the linear function on $\spp ({\R }^{2n+2}, J)$ 
which sends $Z$ to $\trace (YZ)$. The map 
\begin{displaymath}
({\R }^{2n+2}, Y, e_{2n+1};J) \mapsto 
{\ell }_Y|{\spp ({\R }^{2n+2}, J)}_{e_{2n+1}} 
\end{displaymath}
induces a bijection between affine cotypes and coadjoint 
orbits of the odd real symplectic group 
${\Spp ({\R }^{2n+2}, J)}_{e_{2n+1}}$ 
on ${\spp ({\R }^{2n+2}, J)}^{\ast }_{e_{2n+1}}$, the dual of 
its Lie algebra. \medskip  

\noindent \textbf{Proof.} Suppose that the tuples 
$({\R }^{2n+2}, Y, e_{2n+1}; J)$ and $({\R }^{2n+2}, Y', e_{2n+1}; J)$ 
are equivalent. Then there is a real linear map 
$P \in {\Spp ({\R }^{2n+2}, J)}_{e_{2n+1}}$ and a vector $w \in 
{\R }^{2n+2}$ such that $Y' = P(Y + L_{w,e_{2n+1}})P^{-1}$. We 
begin our argument with three observations. \medskip 

\noindent \textbf{Observation 1.} Let $w = w_0 e_0 + \widetilde{w} + w_{2n+1} e_{2n+1} 
\in {\R}^{2n+2}$. Then the matrix of $L_{w, e_{2n+1}}$ 
with respect to the standard basis 
$\mathfrak{e}$ of ${\R }^{2n+2}$ is {\tiny $ \left( 
\begin{array}{ccr} w_0 & 0 & 0 \\ \widetilde{w} & 0 & 0 \\
2w_{2n+1} & {\widetilde{w}}^{\ast } & -w_0 \end{array} 
\right) $}. To see this we compute 
\begin{align}
L_{w,e_{2n+1}}(e_0) & = (w \otimes e^{\ast }_{2n+1})(e_0) + e_{2n+1} \otimes w^{\ast}(e_0) \notag \\
& =  (e^T_0Je_{2n+1})w + (e^T_0Jw)e_{2n+1} 
\notag \\
& = w + w_{2n+1}e_{2n+1} = w_0e_0 + \widetilde{w} + 2w_{2n+1}e_{2n+1}, 
\notag \\
L_{w,e_{2n+1}}(e_i) & = (e^T_i J e_{2n+1}) w + (e^T_i Jw) e_{2n+1} =  {\widetilde{w}}^{\ast }(e_i)e_{2n+1}, \notag \\
L_{w,e_{2n+1}}(e_{2n+1}) & = (e^T_{2n+1}Je_{2n+1})w + (e^T_{2n+1}Jw)e_{2n+1} 
= -w_0 e_{2n+1}.  \notag
\end{align} 
Let $Y =${\tiny $\begin{pmatrix} 
w_0 & -{\widetilde{v}}^{\, \ast } & c \\
\widetilde{w} & \widetilde{Y} & \widetilde{v} \\
2w_{2n+1} & {\widetilde{w}}^{\, \ast } & -w_0 \end{pmatrix}$}$ \in \spp ({\R }^{2n+2}, J)$. Then 
using observation 1 we have
$Y = Y' +L_{w, e_{2n+1}}$, where $Y'=${\tiny $\begin{pmatrix} 
0 & -{\widetilde{v}}^{\, \ast } & c \\
0 & \widetilde{Y} & \widetilde{v} \\
0 & 0 & 0 \end{pmatrix}$} and $w = w_0 e_0 + \widetilde{w} + w_{2n+1}e_{2n+1}$. 
Thus the tuples $({\R }^{2n+2}, Y, e_{2n+1}; J)$ and $({\R }^{2n+2}, Y', e_{2n+1}; J)$ are equivalent. \medskip 

\noindent \textbf{Observation 2.} For $P \in \Spp ({\R }^{2n+2}, J)$ and 
$Y \in \spp ({\R }^{2n+2}, J)$ we have ${\ell }_{PYP^{-1}} = 
{\Ad }^T_{P^{-1}}{\ell }_Y$. To see this, we compute. 
Let $Z \in \spp ({\R }^{2n+2}, J)$. Then  
\begin{align}
{\ell }_{PYP^{-1}}(Z) & = \trace ((PYP^{-1})Z) 
= \trace (P(YP^{-1}Z)) \notag \\
& = \trace ((YP^{-1}Z)P) = \trace (Y(P^{-1}ZP)) = 
{\ell }_Y(P^{-1}ZP) \notag \\
& = {\ell }_Y({\Ad }_{P^{-1}}Z) = ({\Ad }^T_{P^{-1}}{\ell }_Y)Z. \notag
\end{align}

\noindent \textbf{Observation 3.} Let ${\spp ({\R }^{2n+2}, J)}^{\circ}_{e_{2n+1}}$ be the set of 
all $X \in \spp ({\R }^{2n+2},J)$ such that ${\ell }_X(Y) =0$ for 
every $Y \in {\spp ({\R }^{2n+2},J)}_{e_{2n+1}}$. Then
\begin{displaymath}
{\spp ({\R }^{2n+2}, J)}^{\circ}_{e_{2n+1}} = \{ {\ell }_{L_{v,e_{2n+1}}} \in {\spp ({\R }^{2n+2}, J)}^{\ast } 
\setrule \, v \in {\R }^{2n+2} \} .
\end{displaymath}
To see this let $X=${\tiny $\begin{pmatrix} a & -{\widetilde{v}}^{\ast } & c \\
\widetilde{d} & \widetilde{X} & \widetilde{v} \\
2f & {\widetilde{d}}^{\ast } & -a \end{pmatrix}$}$\in \spp ({\R }^{2n+2},J)^{\circ}_{e_{2n+1}}$ and let 
$Y=${\tiny $\begin{pmatrix} 0 & 0 & 0 \\ \widetilde{e} & \widetilde{Y} & 0 \\ 2g & {\widetilde{e}}^{\ast } & 0 
\end{pmatrix}$}$\in {\spp ({\R }^{2n+2}, J)}_{e_{2n+1}}$. Then by hypothesis we have 
\begin{align}
0 & = \mathrm{tr}\mbox{\footnotesize $\begin{pmatrix} a & -{\widetilde{v}}^{\ast } & c \\
\widetilde{d} & \widetilde{X} & \widetilde{v} \\
2f & {\widetilde{d}}^{\ast } & -a \end{pmatrix} \, 
\begin{pmatrix} 0 & 0 & 0 \\ \widetilde{e} & \widetilde{Y} & 0 \\ 2g & {\widetilde{e}}^{\ast } & 0 
\end{pmatrix}$} \notag  \\ 
& = 2gc -{\widetilde{v}}^{\ast }(\widetilde{e}) + \mathrm{tr}(\widetilde{v} \otimes {\widetilde{e}}^{\ast }) 
+ \mathrm{tr}\, \widetilde{X} \widetilde{Y} , \label{eq-onenw} 
\end{align}
for all $g \in \R $, $\widetilde{e} \in {\R }^{2n}$, and $\widetilde{Y} \in \spp ({\R }^{2n}, \widetilde{J})$. Pick 
$g = c$, $\widetilde{e} = 0$, and $\widetilde{Y}=0$. Then (\ref{eq-onenw}) reads $0 = 2c^2$, which 
implies $c=0$. So now (\ref{eq-onenw}) reads
\begin{equation}
0 = -{\widetilde{v}}^{\ast }(\widetilde{e}) + \mathrm{tr}(\widetilde{v} \otimes {\widetilde{e}}^{\ast }) 
+ \mathrm{tr}\, \widetilde{X} \widetilde{Y}, 
\label{eq-twonw}
\end{equation}
for every $\widetilde{Y} \in \spp ({\R }^{2n}, \widetilde{J})$ and every $\widetilde{e} \in {\R }^{2n}$. Pick 
$\widetilde{Y} ={\widetilde{X}}^T \in \spp ({\R }^{2n}, \widetilde{J} )$ and $\widetilde{e} =0$. Then 
(\ref{eq-twonw}) reads $0 = \mathrm{tr}\, \widetilde{X}{\widetilde{X}}^T$, which implies $\widetilde{X} =0$. Now equation (\ref{eq-twonw}) becomes
\begin{equation} 
0 = -{\widetilde{v}}^{\ast }(\widetilde{e}) + \mathrm{tr}( \widetilde{v} \otimes {\widetilde{e}}^{\ast }) = 
{\widetilde{v}}^{\, T} \widetilde{J}\, \widetilde{e} - {\widetilde{e}}^{\, T} \widetilde{J}\, \widetilde{v} = 
2 {\widetilde{v}}^{\, T}\widetilde{J}\, \widetilde{e} = 2 \, \widetilde{\omega }(\widetilde{e}, \widetilde{v} ),  
\label{eq-threenw}
\end{equation}
for every $\widetilde{e} \in {\R }^{2n}$. But $\widetilde{\omega }$ is nondegenerate. Therefore $\widetilde{v} = 
0$. Consequently, ${\spp ({\R }^{2n+2}, J)}^{\circ}_{e_{2n+1}}$ equals 
\begin{align}
&\{ {\ell }_{L_{v, e_{2n+1}}}  \in {\spp ({\R }^{2n+2}, J)}^{\ast } \setrule \, 
v = a e_0 + \widetilde{d} + f e_{2n+1} \in {\R }^{2n+2} \} . \tag*{$\square $}
\end{align} 

\noindent \textbf{Proof of proposition 4.} Suppose that the tuples $({\R }^{2n+2},Y, e_{2n+1}; J)$ and 
$({\R }^{2n+2}, Y', e_{2n+1}; J)$ are equivalent. Then there is $P\in {\Spp ({\R }^{2n+2},J)}_{e_{2n+1}}$ 
and a vector $w \in {\R }^{2n+2}$ such that $Y' = P(Y +L_{w, e_{2n+1}})P^{-1}$. For every 
$Z \in {\spp ({\R }^{2n+2},J)}_{e_{2n+1}}$ we have 
\clearpage 
\begin{align}
{\ell }_{P(Y + L_{w, e_{2n+1}})P^{-1}}(Z) & = 
{\ell }_{PYP^{-1}}(Z) + {\ell }_{PL_{w,e_{2n+1}}P^{-1}}(Z) \notag \\
& \hspace{-.75in} = {\ell }_{PYP^{-1}}(Z) + {\ell }_{L_{Pw,e_{2n+1}}}(Z), \quad 
\mbox{since $P \in {\Spp ({\R }^{2n+2}, J)}_{e_{2n+1}}$} \notag \\
& \hspace{-.75in} = {\ell }_{PYP^{-1}}(Z), \quad \mbox{since $Z \in 
{\spp ({\R }^{2n+2},J)}_{e_{2n+1}}$} \notag \\
& \hspace{-.75in} = ({\Ad }^T_{P^{-1}}{\ell }_Y)Z = ({\Ad }^T_{P^{-1}}({\ell }_Y|
{\spp ({\R }^{2n+2},J)}_{e_{2n+1}}))Z. \notag
\end{align}
Thus the affine cotype, represented by $({\R }^{2n+2}, Y , e_{2n+1};J)$, corresponds to the unique coadjoint orbit of ${\Spp ({\R }^{2n+2}, J)}_{e_{2n+1}}$ through 
${\ell }_Y|{\spp ({\R }^{2n+2},J)}_{e_{2n+1}}$ in 
${\spp ({\R }^{2n+2},J)}^{\ast }_{e_{2n+1}}$. \medskip 

Suppose that for some $Y$, $Y' \in \spp ({\R }^{2n+2}, J)$ and some $P \in {\Spp ({\R }^{2n+2}, J)}_{e_{2n+1}}$
we have ${\ell }_{Y - {\Ad }^T_{P^{-1}}(Y')} =0$ on ${\spp ({\R }^{2n+2},J)}_{e_{2n+1}}$. 
In other words, suppose that ${\ell }_Y|{\spp ({\R }^{2n+2},J)}_{e_{2n+1}}$ lies in the 
${\Spp ({\R }^{2n+2}, J)}_{e_{2n+1}}$ coadjoint orbit \linebreak 
through ${\ell }_Y'|{\spp ({\R }^{2n+2}, J)}_{e_{2n+1}}$. 
Then ${\ell }_{Y - P(Y')P^{-1}} \in {\spp ({\R }^{2n+2},J)}^{\circ}_{e_{2n+1}}$. Therefore 
for some $v \in {\R}^{2+2}$ we have ${\ell }_{Y - P(Y')P^{-1}} = \ell _{L_{v, e_{2n+1}}}$. So 
\begin{displaymath}
Y = P(Y')P^{-1} +L_{v, e_{2n+1}} = P(Y'+ L_{P^{-1}v, e_{2n+1}})P^{-1}. 
\end{displaymath}
Hence the tuples $({\R }^{2n+2}, Y, e_{2n+1}; J)$ and $({\R }^{2n+2}, Y', e_{2n+1}; J)$ are 
equivalent. Thus a coadjoint orbit uniquely determines an affine cotype. Since every 
element of ${\spp ({\R }^{2n+2}, J)}^{\ast }_{e_{2n+1}}$ is of the form 
${\ell }_Y|{\spp ({\R }^{2n+2}, J)}_{e_{2n+1}}$ for some $Y \in \spp ({\R }^{2n+2}, J)$, the map 
from the affine cotype, represented by the tuple $({\R }^{2n+2}, Y, e_{2n+1}; J)$ to 
an ${\Spp ({\R }^{2n+2},J)}_{e_{2n+1}}$ coadjoint orbit through ${\ell }_Y|{\spp ({\R }^{2n+2}, J)}_{e_{2n+1}}$ 
is surjective. This proves proposition 4. \hfill $\square $ \medskip 

Let $(V,Y,v; \omega )$ be a tuple which represents the cotype 
$\nabla $. If $V = V_1 \oplus V_2$, where $v \in V_1$, $V_2 
\ne \{ 0 \}$, and $V_i$ are $Y$-invariant, $\omega $-perpendicular, 
$\omega $ nondegenerate subspaces of $(V,\omega )$, 
then $\nabla $ is the \emph{sum} of 
a cotype $\widetilde{\nabla}$, represented by the tuple 
$(V_1, Y|V_1, v; \omega |V_1)$, and a type $\Delta $, represented 
by the pair $(V_2, Y|V_2; \omega |V_2)$, see \cite{burgoyne-cushman}.  We write 
$\nabla = \widetilde{\nabla} + \Delta$. We say that 
the cotype $\nabla $ is \emph{indecomposable} if it cannot be 
written as the sum of a cotype and a type. If $V_1 = \{ 0 \}$, then 
$v =0$ and the tuple $(\{ 0 \} , 0, 0; 0)$ represents the 
\emph{zero cotype}, which we denote by $\mathbf{0}$. Note that 
$\dim \mathbf{0} =0$. \medskip 

\noindent \textbf{Proposition 5.} An indecomposable affine cotype $\nabla $ is represented by the tuple 
$(V,Y,v^{\circ}; \omega )$, where $Y$ is nilpotent. \medskip 

\noindent \textbf{Proof.} Let $S$ be the semisimple summand in 
the Jordan decomposition of $Y$. The following argument shows that 
$v^0 \in \ker S$. Since $V = \ker S \oplus \mathrm{im}\, S$, 
we may write $v^0 = v_1 + v_2 \in \ker S \oplus \mathrm{im}\, S$. So $0 = Yv^{\circ} = 
Yv_1 + Yv_2 \in \ker S \oplus \mathrm{im}\, S$. Thus $Yv_2 =0$. Since $Y$ is invertible on 
$\mathrm{im}\, S$, it follows that $v_2 =0$. Thus $v^{\circ} = v_1 \in \ker S$. Let $N$ be the nilpotent summand 
of the Jordan decomposition of $Y$. By results of \cite{burgoyne-cushman} the type 
${\Delta }'$, represented by the tuple $(V' = \ker S, N = Y|\ker S; \omega |\ker S)$ is 
a sum of indecomposable types, which is unique up to reordering of the summands. 
Since $v^{\circ} \in \ker N$, there 
is exactly one indecomposable summand ${\Delta}''$, represented by the tuple 
$(V'', N|V''; \omega |V'' )$ where $v^{\circ} \in V''$. So $V''$ is a $Y$ invariant, $\omega $ 
nondegenerate subspace of $(V', \omega |V')$, and hence of $(V, \omega )$, which contains 
$v^{\circ}$. On $V''$ the linear map $Y$ is nilpotent. Thus $(V'', Y|V'', v^{\circ}; {\omega }|V'')$ 
is a tuple, which represents a nilpotent affine cotype ${\nabla }''$. Since the affine cotype $\nabla $ is indecomposable, it follows that the tuples $(V,Y, v^{\circ}; \omega )$ and $(V'', Y''|V'', v^{\circ}; {\omega }|V'')$ are equal. Thus the affine cotype $\nabla $ is nilpotent. \hfill $\square $

\section{Classification of affine cotypes}
%%%%%%%%%%%%%%%%%%%%%%

In this section we classify indecomposable affine cotypes. \medskip 

For $w = w_0 e_0 + \widetilde{w} + w_{2n+1} e_{2n+1} \in 
{\R }^{2n+2}$ the matrix of $L_{w,e_{2n+1}}$ with respect to 
the standard basis $\mathfrak{e}$ is 
{\tiny $ \left( \begin{array}{ccr} w_0 & 0 & 0 \\ \widetilde{w} & 0 & 0 \\
2w_{2n+1} & {\widetilde{w}}^{\ast } & -w_0 \end{array} \right) $}. With 
respect to the basis $\mathfrak{e}$ the matrix of $Y \in  
\spp ({\R }^{2n+2}, J)$ is 
\begin{equation}
\left( \begin{array}{ccc} 
a & -{\widetilde{v}}^{\ast } & c \\
\widetilde{d} & \widetilde{Y} & \widetilde{v} \\
2f & {\widetilde{d}}^{\ast } & -a \end{array} \right) , 
\label{eq-coadtwo}
\end{equation}
where {\footnotesize $\begin{pmatrix} a & c \\ 2f & -a 
\end{pmatrix}$}$\in \spp \big( {\R }^2, ${\tiny 
$\begin{pmatrix} 0 & 1 \\ -1 & 0 \end{pmatrix} $}$\big) $, using the basis $\{ e_0 , e_{2n+1} \}$ of ${\R }^2$. 
Also $\widetilde{d} ,\widetilde{v} \in {\R }^{2n}$. If $({\R }^{2n+2}, 
Y, e_{2n+1}; J)$ is a tuple with the matrix of $Y$ given 
by (\ref{eq-coadtwo}), then we call the ${(0,2n+1)}^{\mathrm{th}}$ 
entry $c$ of $Y$, the \emph{parameter of the tuple}. \medskip 
 
\noindent \textbf{Lemma 6.} Let $({\R }^{2n+2}, Y, e_{2n+1}; J)$ be 
a tuple with parameter $c$, which represents the affine 
cotype $\nabla $. Then $c$ does not depend on the choice of 
representative of $\nabla $. \medskip 

\noindent \textbf{Proof.} Let $({\R }^{2n+2}, Y', e_{2n+1}; J)$ 
be another representative of the cotype $\nabla $. Then 
the tuples $({\R }^{2n+2}, Y, e_{2n+1}; J)$ and 
$({\R }^{2n+2}, Y', e_{2n+1}; J)$ are equivalent. In other words, 
there is a $P \in {\Spp ({\R }^{2n+2}, J)}_{e_{2n+1}}$ and a 
vector $w \in {\R }^{2n+2}$ such that 
\begin{equation}
Y' = P(Y +L_{w,e_{2n+1}})P^{-1}.
\label{eq-coadthree}
\end{equation}
We need only calculate the right hand side of equation 
(\ref{eq-coadthree}). Observe that 
\begin{displaymath}
{\Spp ({\R }^{2n+2}, J)}_{e_{2n+1}} = \Big\{  \mbox{{\footnotesize 
$\left( \begin{array}{ccc} 1 & 0 & 0 \\ \widetilde{u} & 
\widetilde{P} & 0 \\ k & {\widetilde{u}}^{\ast }P & 1 
\end{array} \right) $}} \setrule \, \mbox{{\footnotesize 
$ k \in \R , \widetilde{u} \in {\R }^{2n},  \widetilde{P} \in 
\Spp ({\R}^{2n}, \widetilde{J}) $}} \Big\} . 
\end{displaymath}
Then 
\begin{align}
Y' & = \mbox{{\tiny  $\left( \begin{array}{ccc} 1 & 0 & 0 \\
\widetilde{u} & \widetilde{P} & 0 \\ k & 
{\widetilde{u}}^{\ast }\widetilde{P} & 1 \end{array} \right) $}} \, 
\mbox{{\tiny  $\left( \begin{array}{ccc} a+w_0 & 
-{\widetilde{v}}^{\ast } & c \\ \widetilde{d} + \widetilde{w} &
\widetilde{Y} & \widetilde{v} \\ 2f+2w_{2n+1} & 
(\widetilde{d} +\widetilde{w})^{\ast } & -(a+w_0) \end{array} \right) $}} 
\, \mbox{{\tiny  $\left( \begin{array}{ccc} 1 & 0 & 0 \\
-{\widetilde{P}}^{-1}\widetilde{u} & {\widetilde{P}}^{-1} & 0 \\
-k & -{\widetilde{u}}^{\ast } & 1 \end{array} \right) $}} 
\label{eq-coadfour} \\
& = \mbox{{\footnotesize  $\left( \begin{array}{ccc} * & * &  c \\
* & * & * \\ * & * & * \end{array} \right)$}} \, 
\mbox{{\footnotesize  $\left( \begin{array}{ccc} 1 & 0 & 0 \\
* & * & 0 \\ * & * & 1 \end{array} \right) $}} 
 = \mbox{{\footnotesize  $\left( \begin{array}{ccc} * & * &  c \\
* & * & * \\ * & * & * \end{array} \right)$}.} \notag 
\end{align}
Thus the tuple $({\R }^{2n+2}, Y', e_{2n+1}; J)$ has parameter $c$ 
as well. \hfill $\square $ \medskip 
  
From the conclusion of lemma 6 it makes sense to call $c$ the \emph{parameter of the affine cotype} 
$\nabla $, represented by the tuple $({\R }^{2n+2}, Y, e_{2n+1}; J)$. \medskip 
 
Calculating the right hand side of equation (\ref{eq-coadfour}) 
more explicitly gives
\begin{equation}
Y' = \mbox{{\footnotesize $\left( \begin{array}{ccc}
a' & -(c\widetilde{u} + \widetilde{P}\widetilde{v})^{\ast } & c \\
\rowspace {\widetilde{d}}' & \widetilde{P}\widetilde{Y}{\widetilde{P}}^{-1} 
-L_{\widetilde{u}, \widetilde{P}\widetilde{v}} 
-c\widetilde{u} \otimes {\widetilde{u}}^{\ast } & 
c\widetilde{u} + \widetilde{P}\widetilde{v} \\ 
\rowspace 2f' & ({\widetilde{d}}')^{\ast } & -a' \end{array} \right) $},}
\label{eq-coadfive}
\end{equation}
for some $a'$, $f' \in \R$, ${\widetilde{d}}' \in {\R }^{2n}$, which depend on $k$, $\widetilde{u} \in {\R }^{2n}$ 
and $\widetilde{P}\in \Spp ({\R }^{2n}, \widetilde{J})$. \medskip 

Suppose that the tuple $({\R }^{2n+2},Y, e_{2n+1}; J)$ has parameter 
$c$, which is \emph{nonzero}. Set $k=0$, $\widetilde{P} = I_{2n}$, and $\widetilde{u} = 
-\frac{1}{c}\widetilde{v}$. This determines
matrix of $P \in {\Spp ({\R }^{2n+2}, J)}_{e_{2n+1}}$ in 
equation (\ref{eq-coadfour}). Then equation (\ref{eq-coadfive}) 
reads 
\begin{displaymath}
Y'  = \mbox{\footnotesize $ \left( \begin{array}{ccc} 
a'' & 0 & c \\ 
\rule{0pt}{11pt} {\widetilde{d}}'' & \widetilde{Y} + 
\frac{1}{2c} L_{\widetilde{v}, \widetilde{v}} & 0 \\
\rule{0pt}{11pt} 2f'' & ({\widetilde{d}}'')^{\ast } & - a'' \end{array} \right) $}  
= Y'' + L_{w'',e_{2n+1}}.  
\end{displaymath}
Here $Y''=${\tiny $\begin{pmatrix} 0 & 0 & c \\ 0 & {\widetilde{Y}}' & 0 \\ 0 & 0 & 0 \end{pmatrix}$} with  
${\widetilde{Y}}' = \widetilde{Y} + \frac{1}{2c}L_{\widetilde{v}, \widetilde{v}}$, and 
$w'' = a'' e_0 +{\widetilde{d}}'' + f'' e_{2n+2}$. Therefore the 
tuple $({\R }^{2n+2}, Y', e_{2n+1};J)$ is equivalent to the 
tuple $({\R }^{2n+2}, Y'', $ $e_{2n+1};J)$. Thus we have proved \medskip 

\noindent \textbf{Proposition 7.} If the affine cotype $\nabla $, represented by the tuple 
$({\R}^{2n+2},Y',$ $e_{2n+1}; J)$, has 
a nonzero parameter $c$, then it is decomposable into a 
sum of a two dimensional indecomposable nilpotent affine cotype ${\nabla }_1(0), c$, 
represented by the tuple $\big({\R }^2 $, \raisebox{1pt}{\tiny $\begin{pmatrix} 0 & c \\ 0 & 0 \end{pmatrix} $}, $e_{2n+1}$; \raisebox{1pt}{\tiny $\begin{pmatrix} 0 & 1 \\ -1 & 0 \end{pmatrix} $}$\big)$ 
with basis $\{  e_0, e_{2n+1} \} $ of ${\R }^2$, and a type $\Delta $, represented by the pair $({\R }^{2n}, {\widetilde{Y}}'; \widetilde{J})$.  \medskip 

We can sharpen the conclusion of proposition 7 a bit. By proposition 4 of \cite{burgoyne-cushman} the height of the type $\Delta $ in proposition 7 above is strictly less than the height of ${\Delta }_1, c$, which is $1$. Thus the height of $\Delta $ is $0$, that is, $\Delta $ is a semisimple type. \medskip 

Suppose that the affine cotype $\nabla $, represented by the tuple $({\R }^{2n+2}, Y,$ \linebreak 
$e_{2n+1};J)$, where $Y$ is given in equation (\ref{eq-coadtwo}), has parameter equal 
to $0$. The tuple $({\R }^{2n}, \widetilde{Y}, \widetilde{v}; \widetilde{J})$ represents the \emph{little cotype} ${\nabla }_{\ell }$ associated to $\nabla $. \medskip 

\noindent \textbf{Lemma 8.} The little cotype ${\nabla }_{\ell }$ is uniquely 
determined by the affine cotype $\nabla $ with parameter $0$. \medskip 

\noindent \textbf{Proof}. Let $({\R }^{2n+2}, Y', e_{2n+1}; J)$ 
be another representation of the affine cotype $\nabla $ with parameter $0$. Here  
\begin{equation}
Y' = \mbox{{\footnotesize $\left( \begin{array}{ccc} 
a' & -(\widetilde{P}\widetilde{v})^{\ast } & 0 \\
{\widetilde{d}}' & \widetilde{P}\widetilde{Y}{\widetilde{P}}^{-1} 
-L_{\widetilde{u}, \widetilde{P}\widetilde{v}} & 
\widetilde{P}\widetilde{v} \\
2f' & ({\widetilde{d}}')^{\ast } & -a' \end{array} \right) $},}
\label{eq-coadsix}
\end{equation}
for some $\widetilde{u} \in {\R }^{2n}$. Equation (\ref{eq-coadsix}) is obtained from (\ref{eq-coadfive}) by setting $c=0$. The tuple $({\R }^{2n}, {\widetilde{Y}}',\widetilde{P}\widetilde{v}; \widetilde{J})$, where 
${\widetilde{Y}}' = \widetilde{P}(\widetilde{Y} - L_{{\widetilde{P}}^{-1}\widetilde{u}, \widetilde{v}}) 
{\widetilde{P}}^{-1}$, is equivalent to the tuple 
$({\R }^{2n}, \widetilde{Y}, \widetilde{v}; \widetilde{J})$. Hence it also represents ${\nabla }_{\ell }$. 
\hfill $\square $ \medskip 

\noindent \textbf{Lemma 9.} Let $\nabla $ be an affine cotype with parameter 
equal to $0$. Then $\nabla $ is uniquely determined by its 
little cotype ${\nabla }_{\ell }$. \medskip 

\noindent \textbf{Proof}. Suppose that the affine cotypes 
$\nabla $ and ${\nabla }'$ with parameter zero, \linebreak 
represented by the tuples $({\R}^{2n+2},Y,e_{2n+1};J)$ and 
$({\R}^{2n+2},Y',e_{2n+1};J)$, \linebreak  
respectively, both have the same little cotype ${\nabla }_{\ell}$, 
represented by $({\R }^{2n},\widetilde{Y}, \widetilde{v}; 
\widetilde{J})$. Say that 
\begin{displaymath}
Y= \mbox{{\footnotesize $\left( \begin{array}{crc} 
w_0 & -{\widetilde{v}}^{\ast } & 0 \\
\widetilde{w} & \widetilde{Y} & \widetilde{v} \\
2w_{2n+1} & {\widetilde{w}}^{\ast } & -w_0 \end{array} \right) $}}
\, \, \, \mathrm{and} \, \, \, 
Y'= \mbox{{\footnotesize $\left( \begin{array}{crc} 
w'_0 &-(\widetilde{v}')^{\ast } & 0 \\
\widetilde{w}' & \widetilde{Y}' & \widetilde{v}' \\
2w'_{2n+1} & ({\widetilde{w}}')^{\ast } & -w'_0 \end{array} \right) $}},
\end{displaymath}
where $w = w_0 e_0 + \widetilde{w} + w_{2n+1}e_{2n+1}$ and 
$w' = w'_0 e_0 + \widetilde{w}' + w'_{2n+1}e_{2n+1}$ lie in 
${\R }^{2n}$ and $\widetilde{Y}, \widetilde{Y}'$ lie in 
$\spp ({\R }^{2n}, \widetilde{J})$. By hypothesis the tuples 
$({\R }^{2n}, \widetilde{Y}, \widetilde{v}; \widetilde{J})$ and 
$({\R }^{2n}, \widetilde{Y}', \widetilde{v}'; \widetilde{J})$ are 
equivalent. In other words, there is $\widetilde{P} \in 
\Spp ({\R }^{2n}, \widetilde{J})$ such that 
$\widetilde{P}\widetilde{v} = \widetilde{v}'$ and a vector $\widetilde{w} \in {\R}^{2n}$ such that $\widetilde{Y}' = 
\widetilde{P}(\widetilde{Y} +L_{\widetilde{w}, \widetilde{v}}){\widetilde{P}}^{-1}$. Let $P =
${\tiny $\left( \begin{array}{ccc} 1 & 0 & 0 \\ 
\widetilde{y} & \widetilde{P} & 0 \\ 0 & {\widetilde{y}}^{\, \ast }\widetilde{P} & 1 \end{array} \right) $}. 
Then $P \in {\Spp ({\R }^{2n+2}, J)}_{e_{2n+1}}$. With $c=0$, 
$k=0$, and $\widetilde{u} = \widetilde{y}$ equation (\ref{eq-coadfive}) reads 
\begin{align}
PYP^{-1} & = 
\mbox{{\footnotesize $ \left( \begin{array}{ccc} 
{\widetilde{a}}^{\dagger} & -(\widetilde{P}\widetilde{v})^{\ast} & 0 \\ 
{\widetilde{d}}^{\, \dagger} & \widetilde{P}\widetilde{Y}{\widetilde{P}}^{-1} 
-L_{\widetilde{y}, \widetilde{P}\widetilde{v}} & \widetilde{P}\widetilde{v} \\ 
2{\widetilde{f}}^{\, \dagger} & ({\widetilde{d}}^{\, \dagger})^{\ast } & -{\widetilde{a}}^{\dagger} \end{array} \right) $}} \, 
\label{eq-coadsixstar}
\end{align}
Setting $\widetilde{y} = -\widetilde{P}\widetilde{w}$, 
$Pw = {\widetilde{a}}^{\dagger } e_0 + {\widetilde{d}}^{\, \dagger} + 2{\widetilde{f}}^{\, \dagger} e_{2n+1}$, 
and using the fact that ${\widetilde{Y}}' = \widetilde{P}(Y + L_{\widetilde{w}, \widetilde{v}}){\widetilde{P}}^{-1}$ by hypothesis, equation (\ref{eq-coadsixstar}) reads 
\begin{align}
PYP^{-1} & = \mbox{\footnotesize $\begin{pmatrix} 
0 & -({\widetilde{v}}')^{\ast } & 0 \\
0 & {\widetilde{Y}}' & {\widetilde{v}}' \\
0 & 0 & 0 \end{pmatrix}$}  +PL_{w, e_{2n+1}}P^{-1}  
\label{eq-coadsevenstar} \\
& = \mbox{\footnotesize $\begin{pmatrix}
w'_0 & -({\widetilde{v}}')^{\ast } & 0 \\
{\widetilde{w}}' & {\widetilde{Y}}' & {\widetilde{v}}' \\
2w'_{2n+1} & -({\widetilde{w}}')^{\ast } & -w'_0 \end{pmatrix} $} - L_{w', e_{2n+1}} + L_{Pw, e_{2n+1}} 
\notag \\
& = Y' - P(L_{P^{-1}w'-w, e_{2n+1} })P^{-1}. \notag
\end{align} 
Therefore the tuples $({\R }^{2n+2}, Y, e_{2n+1};J)$ and 
$({\R }^{2n+2}, Y', e_{2n+1};J)$ are equivalent. So the 
affine cotypes $\nabla $ and ${\nabla }'$ with parameter $0$ 
are equal. \hfill $\square $ \medskip 

\vspace{-.15in}\noindent \textbf{Lemma 10.} If $\nabla $ is \emph{not} an affine cotype, 
then $\nabla $ is the sum of the zero cotype $\mathbf{0}$ and a type. \medskip 

\noindent \textbf{Proof}. Suppose that $\nabla $ is represented 
by the tuple $(V, Y, v; \omega )$. Since $\nabla $ is not 
affine, $v=0$. Because $\{ 0 \} $ and $V$ are $Y$-invariant, 
$\omega $-perpendicular, $\omega $-symplectic subspaces 
of $(V, \omega )$ with 
$V = \{ 0 \} \oplus V$ and $0 \in \{ 0 \}$, we can write 
the tuple $(V,Y,0; \omega )$ as the sum of the tuple 
$(\{ 0 \} , 0, 0; 0)$ and the pair $(V,Y; \omega )$. Therefore 
$\nabla $ is the sum of the zero cotype and a type, represented by 
the pair $(V,Y; \omega )$.\mbox{}\hfill $\square $ \medskip 

\noindent \textbf{Proposition 11.} Every affine cotype is either the sum of a 
nonzero nilpotent indecomposable affine cotype and a sum of indecomposable types or 
the sum of the zero cotype and a sum of indecomposable types. This 
decomposition is unique up to reordering of the summands which 
are types. \medskip

\vspace{-.15in}\noindent \textbf{Proof}. Let $\nabla $ be an affine cotype. 
Then by lemma 6 it has a parameter, say $c$. If $c\ne 0$, 
then by proposition 7 the cotype $\nabla $ is the sum 
of the two dimensional nilpotent indecomposable affine cotype ${\nabla }_1(0), c$ 
and a semisimple type $\Delta $, which may be decomposed into a sum of indecomposable semisimple types that are unique up to reordering of the summands, see \cite[table 3, p.349]{burgoyne-cushman}. The argument stops here. 
If $c=0$, then by lemma 8 the affine 
cotype $\nabla $ has a unique little cotype ${\nabla }_{\ell }$, 
which uniquely determines $\nabla $. This correspondence 
respects decomposition, that is, if $\nabla $ is decomposable, 
then so is ${\nabla }_{\ell }$ and conversely. If ${\nabla }_{\ell }$ is 
not affine, then by lemma 10 it is the sum of the zero cotype and 
a type $\Delta $, which may be written as a sum of indecomposable 
types. The argument stops here.  Otherwise ${\nabla }_{\ell }$ 
is affine and we may repeat the above argument. Only a finite number of repetitions are needed before the argument comes to a stop, because  
$\dim {\nabla }_{\ell } < \dim \nabla < \infty $. \hfill $\square $  \medskip  

\noindent \textbf{Proposition 12.} Let $\nabla $ be a nilpotent indecomposable affine 
cotype of dimension $2r+2$ with $r \le n$. Then exactly one of the following 
possibilities holds. \medskip 

\noindent 1. Suppose that $r =0$ and the parameter $d$ is 
\emph{not} equal to $0$. Then 
$({\R }^2, Y, e_1;$\raisebox{1pt}{\tiny $\begin{pmatrix}
0 & 1 \\ -1 & 0 \end{pmatrix}$}$)$ is a representative
of $\nabla $. Moreover, there is  
a basis ${\mathfrak{e}}' =\{ e_0, e_1 \}$ of ${\R }^2$ 
such that the matrix of 
$Y$ with respect to the basis ${\mathfrak{e}}'$ is 
{\tiny $ \left(\begin{array}{rc}
0 & d \\ 0 & 0 \end{array} \right) $}. Note $Y^2 =0$. So $Y$ has height $1$. The 
matrix of $\omega $ with respect to the basis ${\mathfrak{e}}'$ is 
{\tiny $\begin{pmatrix} 0 & 1 \\ -1 & 0 \end{pmatrix}$}. 
We denote this nilpotent indecomposable 
affine cotype by ${\nabla }_1(0),d$. Here the $d$ is nonzero and 
is called a \emph{modulus}.  \medskip 

\noindent 2. Suppose that $r =0$ and the parameter $d$ is equal 
to $0$. Then $({\R }^2, Y, e_1;$ \linebreak 
{\tiny $\begin{pmatrix}
0 & 1 \\ -1 & 0 \end{pmatrix}$}$)$ is a representative
of $\nabla $. Moreover, there is  
a basis ${\mathfrak{e}}' = \{ e_0, e_1 \}$ of ${\R }^2$ 
such that the matrix of $Y$ with respect to the basis ${\mathfrak{e}}'$ is 
\raisebox{2pt}{\tiny $\begin{pmatrix}
0 & 0 \\ 0 & 0 \end{pmatrix}$}. The 
matrix of $\omega $ with respect to the basis ${\mathfrak{e}}'$  is 
\raisebox{2pt}{\tiny $\begin{pmatrix} 0 & 1 \\ -1 & 0 \end{pmatrix}$}. 
We denote this nilpotent indecomposable 
affine cotype by ${\Delta }_0(0,0)$. There is \emph{no} modulus and the parameter 
of ${\Delta}_0(0,0)$ is $0$. \medskip 

\noindent 3. $\nabla $ is obtained by taking 
$r$-iterated little cotypes which ends at the cotype 
${\nabla}_1(0), \, d$ with $d\ne 0$. 
Let $({\R }^{2r+2}, Y,$ $e_{2r+2};J)$ 
be a representative of $\nabla $. Then there is a vector 
$z \in {\R }^{2r+2}$ and a basis 
\begin{align}
\mathfrak{f} & = \{ Y^{2r+1}z; (-1)^rY^{r+1}z, (-1)^{r-1}Y^{r+2}z, \ldots , -Y^{2r}z; \notag \\
&\hspace{.5in} (-1)^rdY^rz, (-1)^rdY^{r-1}z, \ldots ; (-1)^rd z \} 
\notag 
\end{align}
of ${\R }^{2r+2}$ such that $Y^{2r+2} =0$. The matrix of $Y$ with 
respect to the basis $\mathfrak{f} $ is 
\begin{displaymath}
\mbox{{\footnotesize $\left( \begin{array}{c|c|c|c}
0 & -e^T_r & 0 & 0 \\ \hline
\rule{0pt}{10pt} 0 & -N^T   & D & 0 \\ \hline 
0 &    0   & N & e_r \\ \hline
0 &    0   & 0 & 0   
\end{array} \right) $},} 
\end{displaymath}
where $N$ is the $r\times r$ upper Jordan block 
\begin{displaymath}
\mbox{{\tiny $\left( \begin{array}{ccccc} 
0 &1 &       &        &     \\
  &0 & 1     &        &     \\
  &  &\ddots & \ddots &      \\
  &  &       &        & 1    \\
  &  &       &        & 0 
\end{array} \right) $}}
\end{displaymath}
and $D$ is the $r\times r$ antidiagonal matrix {\tiny $\begin{pmatrix} & & & 0 \\ & &  & \\ & 0 &  & \\ d &  & &  
\end{pmatrix} $}. $Y$ is nilpotent of height $2r+1$ and 
${\R }^{2r+2}$ is spanned by one Jordan chain of length $2r+2$. 
The matrix of $\omega $ with respect to the basis $\mathfrak{f}$ is 
\begin{displaymath}
J = \mbox{{\footnotesize $\left( \begin{array}{r|r|c|c}
0 &    0   &   0 &  1 \\ \hline
0 &    0   & I_r & 0 \\ \hline 
0 & -I_r   &   0 & 0  \\ \hline
-1 &    0   &   0 & 0   
\end{array} \right) $},} 
\end{displaymath}
where $I_r$ is the $r\times r$ diagonal matrix 
$\mathrm{diag}(1,1,\ldots ,1)$. We denote this nilpotent indecomposable 
affine cotype by ${\nabla}_{2r+1}(0),d$. 
Here $d$ is a nonzero \emph{modulus} and the parameter of ${\nabla}_{2r+1}(0)$ is zero.  \medskip 

\noindent 4. $\nabla $ is obtained by taking 
$r$-iterated little cotype which ends at the 
cotype ${\nabla }_0(0,0)$. Let $({\R }^{2r+2}, Y,$ $e_{2r+2};J)$ be a representative of $\nabla $. Then there are vectors $z,w \in {\R }^{2r+2}$ and a basis 
\begin{align}
{\mathfrak{f}}'& = \{ (-1)^rY^rz; z, -Yz, \ldots ,
(-1)^{r-1}Y^{r-1}z, Y^rw, Y^{r-1}w, \ldots ; w \} 
\notag 
\end{align}
of ${\R }^{2r+2}$ such that $Y^{r+1}z = Y^{r+1}w =0$. The matrix 
of $Y$ with respect to the basis ${\mathfrak{f}}'$ is 
\begin{displaymath}
\mbox{{\footnotesize $\left( \begin{array}{c|c|c|c}
0 & -e^T_r & 0 & 0 \\ \hline
\rule{0pt}{10pt} 0 & -N^T   & 0 & 0 \\ \hline 
0 &    0   & N & e_{2r} \\ \hline
0 &    0   & 0 & 0   
\end{array} \right) $},} 
\end{displaymath}
where $N$ is the $r\times r$ upper Jordan block. 
$Y$ is nilpotent of height $r$ and 
${\R }^{2r+2}$ is spanned by two Jordan chains of length $r+1$. 
The matrix of $\omega $ with respect to the basis ${\mathfrak{f}}'$ is 
$J$. Denote this nilpotent indecomposable affine cotype by 
${\nabla}_{r}(0,0)$. There is \emph{no} modulus and the 
parameter of ${\nabla }_{r}(0,0)$ is $0$. \medskip 

\noindent \textbf{Proof.} 
\par\noindent 1. and 2. The cotypes ${\nabla }_1(0),d$ and 
${\Delta }_0(0,0)$ are clearly indecomposable. \medskip 

Suppose that $r\ge 1$. Then according to the proof of proposition 11, the indecomposable 
affine cotype $\nabla $ is obtained by taking 
the $r$-iterated little cotype which ends in either a) the cotype 
${\nabla}_1(0), \, d$ with $d\ne 0$ or b) the 
cotype ${\Delta }_0(0,0)$. \medskip 

\par \noindent 3. Alternative a) holds. Using the basis $\widetilde{\mathfrak{f}}$ given by
\begin{align*}
&\mbox{\footnotesize $ \{  Y^{2r+1}z, -Y^{2r}z, Y^{2r-1}z, \ldots , 
(-1)^rY^{r+1}z, (-1)^rdY^rz, (-1)^rdY^{r-1}z, \ldots , (-1)^rd z \} $, }
\end{align*}
instead of the basis $\mathfrak{f}$, we see that the matrix of $Y$ and 
$\omega $ are respectively  
\begin{displaymath}
\mbox{{\tiny $\left( \begin{array}{crrr|cccc}
0 & -1  &        &      &  &   &        & \\
  & 0   & -1     &      &  &   &        & \\
  &     & \ddots &  -1  &  &   &        &  \\
  &     &        &   0  & d &   &        & \\ \hline 
  &     &        &      &0 & 1 &        &  \\
  &     &        &      &  & 0 & 1      &   \\
  &     &        &      &  &   & \ddots & 1 \\
  &     &        &      &  &   &        & 0 
\end{array} \right) $}} \, \, \, \mathrm{and} \, \, \,  
\mbox{{\tiny $\left( \begin{array}{crrr|cccc}
    &   &    &     &  &   &   & 1  \\
    &   &    &     &  &   & 1 &    \\
    &   &    &     &  &   &   &    \\
    &   &    &     & 1 &  &   &    \\ \hline 
    &   &    & -1  &  &   &   &    \\
    &   &-1  &     &  &   &   &   \\
    &   &    &     &  &   &   &   \\
-1  &   &    &     &  &   &   &  
\end{array} \right) $}.}
\end{displaymath}
From this it is clear that the little cotype of the affine 
cotype ${\nabla }_{2r+1}(0),d$ is the affine cotype 
${\nabla }_{2r-1}(0),d$. Because ${\nabla}_1(0),d$ is an 
indecomposable cotype, it follows that ${\nabla }_{2r+1}(0),d$ is 
also. Note that $Y$ is nilpotent of height $2r+1$ and that 
${\R }^{2r+2}$ is spanned by one Jordan chain of length $2r+2$. \medskip 

\par \noindent 4. Alternative b) holds. Using the basis $\widetilde{\mathfrak{f}}\, '$ given by
\begin{align*}
& \mbox{\footnotesize $ \{  Y^rz, -Y^{r-1}z, Y^{r-2}z, \ldots , (-1)^rz, 
%& \hspace{.5in}
(-1)^rY^rw, (-1)^rY^{r-1}w, \ldots , (-1)^rw \} $}, 
\end{align*}
instead of the basis ${\mathfrak{f}}\, '$, we see that the matrix of $Y$ and 
$\omega $ are respectively 
\begin{displaymath}
\mbox{{\tiny $\left( \begin{array}{crrr|cccc}
0 & -1  &        &      &  &   &        & \\
  & 0   & -1     &      &  &   &        & \\
  &     & \ddots &  -1  &  &   &        &  \\
  &     &        &   0  &  &   &        & \\ \hline 
  &     &        &      &0 & 1 &        &  \\
  &     &        &      &  & 0 & 1      &   \\
  &     &        &      &  &   & \ddots & 1 \\
  &     &        &      &  &   &        & 0 
\end{array} \right) $}} \, \, \, \mathrm{and} \, \, \,  
\mbox{{\tiny $\left( \begin{array}{crrr|cccc}
    &   &    &     &  &   &   & 1  \\
    &   &    &     &  &   & 1 &    \\
    &   &    &     &  &   &   &    \\
    &   &    &     & 1 &  &   &    \\ \hline 
    &   &    & -1  &  &   &   &    \\
    &   &-1  &     &  &   &   &   \\
    &   &    &     &  &   &   &   \\
-1  &   &    &     &  &   &   &  
\end{array} \right) $}.}
\end{displaymath}
From this it is clear that the little cotype of the affine 
cotype ${\nabla }_r(0,0)$ is the affine cotype 
${\nabla }_{r-1}(0,0)$. Because ${\nabla}_0(0,0)$ is an 
indecomposable cotype, it follows that ${\nabla }_r(0,0)$ is 
also. Note that $Y$ is nilpotent of height $r$ and that 
${\R }^{2r+2}$ is spanned by two Jordan chains both of 
length $r+1$. \hfill $\square $

\end{document}